\documentclass[12pt]{article}
\usepackage{titlesec}

\titleformat*{\section}{\sf\normalsize\bfseries}
\titleformat*{\subsection}{\Large\bfseries}
\titleformat*{\subsubsection}{\large\bfseries}
\titleformat*{\paragraph}{\large\bfseries}
\titleformat*{\subparagraph}{\large\bfseries}

\usepackage[T1]{fontenc}
\usepackage[utf8]{inputenc}
\usepackage{amsmath,amssymb,amsfonts,enumitem,mathtools,calligra,mathrsfs}

\usepackage{wasysym}

\usepackage{upgreek}

\setenumerate[1]{label=(\roman*),leftmargin=*,nolistsep,parsep=\parskip}

\usepackage[noadjust]{cite} % noadjust removes the space in front of a citation.

\usepackage[hidelinks]{hyperref}

\usepackage{verbatim}

\usepackage[thmmarks,hyperref,amsmath]{ntheorem}
\theoremstyle{plain}
\newtheorem{thm}{Theorem}[section]

\theorembodyfont{\upshape}
\newtheorem{dfn}{Definition}[section]
\newtheorem{prop}[dfn]{Proposition}

\newtheorem{corollary}[dfn]{Corollary}
\newtheorem{lemma}[dfn]{Lemma}

\theoremsymbol{\ensuremath\Box}
\theoremheaderfont{\scshape}
\theoremseparator{:}
\newtheorem*{proof}{Proof}

\usepackage{float}
\usepackage{tikz}
\usepackage{tikz-cd}
\usetikzlibrary{calc}
\usetikzlibrary{positioning}
\usetikzlibrary{shapes.geometric}
\usetikzlibrary{arrows}
\usetikzlibrary{decorations.markings}
	\tikzstyle{filledNode}=[shape=circle,scale=0.3, draw, fill=black!100]
	\tikzstyle{emptyNode}=[shape=circle,scale=0.3, draw]
	\tikzstyle{invisibleNode}=[scale=1, draw=none, inner sep = 0, outer sep = 0, fill=none]
	\tikzstyle{borderlessNode}=[scale=1, draw=none, inner sep = 0, outer sep = 0]
	\tikzstyle{directedEdge<}=[decoration={markings,mark=at position 0.5 with {\arrow[scale=1.5, >=stealth]{<}}}]
	\tikzstyle{directedEdge>}=[decoration={markings,mark=at position 0.5 with {\arrow[scale=1.5, >=stealth]{>}}}]

\usepackage[font={small}]{caption}

\usepackage{listings}

\usepackage[boxruled]{algorithm2e}

\newcommand{\CC}{\mathbb{C}}

\newcommand{\ZZ}{\mathbb{Z}}
\newcommand{\NN}{\mathbb{N}}

\newcommand{\PP}{\mathbb{P}}

\newcommand{\tensor}{\otimes}

\newcommand{\intersection}{\cap}

\newcommand{\dual}{\vee}
\newcommand{\Pic}{\mathrm{Pic}}

\newcommand{\Hom}{\mathrm{Hom}}

\newcommand{\interior}[1]{\mathring{\kern0pt#1}}

\author{Patrick Bloß}
\date{\today}
\title{The Infinitesimal Torelli Theorem for hypersurfaces in abelian varieties}

\begin{document}
\pagestyle{plain}

\maketitle

\begin{abstract}
	Given a compact Kähler manifold, the Infinitesimal Torelli problem asks whether the differential of the period map of a Kuranishi family is injective. Unlike the classical Torelli theorem for curves, there is a negative answer for example for hyperelliptic curves of genus greater than $2$. Nevertheless the Infinitesimal Torelli Theorem holds for many other classes of manifolds. We will prove it for smooth hypersurfaces in simple abelian varieties with sufficiently high self-intersection giving an effective bound on a result by Green in this particular case.
		
\end{abstract}

\section{Introduction}

Consider a family of compact Kähler manifolds $\phi\colon \mathcal{X} \to B$, i.e. a proper holomorphic submersion of complex manifolds with Kähler fibres. Denote by $X_b$ the fibre $\phi^{-1}(b)$ of $\phi$ over $b\in B$ and fix $0 \in B$. Write $X=X_0$. Ehresmann's theorem ensures that in some neighborhood $U$ of $0$ there are well defined isomorphisms of the cohomology groups $H^k(X_b,\ZZ) \cong H^k(X_0,\ZZ)$. These will in general not preserve the Hodge structure so it makes sense to consider the \textit{period map}. For given $k$ and $p$ the $p$-th piece of the period map with respect to the $k$-th cohomology group is defined by
\[
	\mathcal{P}^{p,k}\colon U \to \mathrm{Grass}(b^{p,k},H^k(X,\CC)), b \mapsto F^pH^k(X_b,\CC) 
\]
where $F^pH^k(X_b,\CC)$ denotes the $p$-th step of the Hodge filtration and $b^{p,k}=\dim F^pH^k(X_b,\CC)$ (note that this is independent of $b$ as all $X_b$ have the same Hodge numbers). Griffiths showed that this map is holomorphic so we can consider its differential 
\[
	d\mathcal{P}^{p,k}\colon T_{B,0} \to \Hom(F^p H^k(X,\CC),H^k(X,\CC)/F^p H^k(X,\CC)). 
\]
Furthermore he showed that $d\mathcal{P}^{p,k}$ is the composition of the Kodaira-Spencer map $T_{B,0} \to H^1(X,T_X)$ with the map
\[
	H^1(X,T_X) \to \Hom(H^{k-p}(X,\Omega_X^p),H^{k-p+1}(X,\Omega_X^{p-1}))
\]
given by the cup product and the interior product.
Now we say that the Infinitesimal Torelli Theorem (ITT in the following) holds for a compact Kähler manifold $X$ if the period map $\mathcal{P}^{n,n}$ of a Kuranishi family of $X$ is an immersion. 
% Immersion vs differential injective on the deformation space
Since the Kodaira-Spencer map is an isomorphism for a Kuranishi family we need to show injectivity of the map
\[
	H^1(X,T_X) \to \Hom(H^0(X,\omega_X),H^1(X,\Omega_X^{n-1})).
\]
For a curve $C$ it follows easily from a classical result by Noether that the ITT holds if and only if $C$ has genus $g(C)=2$ or $g(C)> 2$ and $C$ is non-hyperelliptic. That is to say that in this case very ampleness of the canonical sheaf is a sufficient condition.
% Reference
For surfaces, however, Garra and Zucconi show that for any $n \geq 5$ there exists a generically smooth $n+9$ dimensional irreducible component of the moduli space of algebraic surfaces such that for a general element of it the ITT fails (see \cite{MR2413341}). Thus finding classes of objects that satisfy the ITT is still an open problem. Reider proves it for surfaces of irregularity at least $5$ with globally generated cotangent bundle that satisfy some additional conditions (see \cite{MR929539}).

The ITT has been shown by Griffiths to hold for hypersurfaces in projective space. 
Green then generalized this to \textit{sufficiently ample} hypersurfaces in an arbitrary smooth projective variety $Y$. By sufficiently ample he means that there exists an ample line bundle $L_0$ such that it holds for all smooth hypersurfaces $X\subset Y$ with $\mathcal{O}_Y(X) \tensor L_0^{-1}$ ample. He does however not give an effective bound on how ample $L_0$ has to be.

We consider specifically the case where $Y=A$ is an abelian variety.
Our main result is the following.
\begin{thm}
\label{thm_main}
	Let $X$ be a smooth hypersurface in a $g$-dimensional simple abelian variety $A$. If 
	$h^0(A,\mathcal{O}_A(X))>\left(\frac{g}{g-1}\right)^g \cdot g!$, then the Infinitesimal Torelli Theorem holds for $X$.
\end{thm}
Following Green's method one can see that for a hypersurface $X \subset A$ in a $g$-dimensional variety, if $L \coloneqq \mathcal{O}_A(X)$ is ample, the ITT holds if the mutliplication map
\[
	H^0(A,L) \tensor H^0(A,L^{g-1}) \to H^0(A,L^g)
\]
is surjective. We will discuss this in more detail in Section \ref{section_green}. Closely related to this is the notion of projective normality. An ample line bundle $L$ on an abelian variety $A$ is very ample and defines a projectively normal embedding if the mutliplication map
\[
H^0(A,L) \tensor H^0(A,L) \to H^0(A,L^2)
\]
is surjective. By an inductive argument it is easy to see that in this case the ITT holds as well. It is well known that $L$ defines a projectively normal embedding if $L=M^n$ with $n \geq 3$ or $n=2$ and some additional condition on the basepoints of $L$ is satisfied. Projective normality for primitive line bundles is fully understood for abelian surfaces (see for example \cite{MR3656291}). For the higher dimensional case Hwang and To give a bound for projective normality to hold for a very general polarized abelian variety in terms of the self-intersection of the line bundle (see \cite{MR2964474}). Finally Iyer gives a bound for higher dimensional simple abelian varieties (see \cite{MR1974682}). In Section \ref{section_main_thm} we will use this approach to prove the following theorem. 
\begin{thm}
	\label{thm_main_2}
	Let $L$ be a line bundle on a simple abelian variety $A$ of dimension $g$.
	Fix $n \in \NN$. If $h^0(A,L)>\left(\frac{n+1}{n}\right)^g \cdot g!$ then the multiplication map
	\[
	\mu_n\colon H^0(A,L) \tensor H^0(A,L^n) \rightarrow H^0(A,L^{n+1})
	\]
	is surjective.
\end{thm}
Theorem \ref{thm_main} is then a corollary of this.

For non-simple abelian threefolds \cite{2018arXiv180308780L} gives numerical conditions for projective normality, taking into account all possible abelian subvarieties.

\section{Green's Approach}
\label{section_green}
Let $X$ be a smooth ample hypersurface in an arbitary smooth projective variety $Y$ of dimension $d$. Let $L \coloneqq \mathcal{O}_Y(X)$. There is a short exact sequence
\[
	0 \to N_X^\dual \to \Omega_Y^1 \tensor \mathcal{O}_X \to \Omega_X^1 \to 0
\]
where $N_X$ denotes the normal bundle of $X$ in $Y$. For any $p\geq1$ this gives a long exact sequence
\[
	0 \to S^p N_X^\dual \to \ldots \to \Omega_Y^{p-1} \tensor N_X^\dual \to \Omega_Y^p \tensor \mathcal{O}_X \to \Omega_X^1 \to 0.
\]
Green then obtains a spectral sequence abutting to zero from which he ultimately deduces (under the assumption that $L$ is sufficiently ample) the following commutative diagram
\begin{equation}
\label{cd_1}
	\begin{tikzcd}
		H^0(X,\omega_X) \arrow{r} \tensor H^1(X,\Omega_X^{d-1})^\dual & 
		H^1(X,T_X)^\dual \\
		H^0(X,\omega_X) \tensor H^0(X,L|_X^{d-1} \tensor \omega_X) \arrow{r} \arrow{u} &
		H^0(X,L|_X^{d-1} \tensor \omega_X^2) \arrow[twoheadrightarrow]{u} \\
		H^0(Y,L \tensor \omega_Y) \tensor H^0(Y,L^d \tensor \omega_Y) \arrow{r} \arrow[twoheadrightarrow]{u} &
		H^0(Y,L^{d+1} \tensor \omega_Y^2) \arrow[twoheadrightarrow]{u}		
	\end{tikzcd}
\end{equation}
The two vertical maps on the bottom are simply restriction maps. Their surjectivity is obtained from the vanishing of certain cohomology groups. The vertical map on the top right comes from a quotient map and is thus surjective as well. Finally the map on the top is the dual of $d\mathcal{P}^{d,d}$ so the ITT holds if the mutliplication map on the bottom is surjective. Note that the ITT may still hold if this map fails to be surjective.

Consider the product $Y \times Y$ and denote by $\pi_i\colon Y \times Y \to Y$ for $i=1,2$ the projection maps to the two factors. Furthermore let $\Delta = \{(y,y) \ | \ y \in Y \}$ be the diagonal in $Y \times Y$ and let $\mathcal{I}_{\Delta/Y}$ denote its ideal sheaf. Under the assumption that $L$ is sufficiently ample Green then deduces the surjectivity of the map on the bottom of diagram (\ref{cd_1}) from the vanishing of 
$H^1(Y \times Y, \mathcal{I}_{\Delta/Y} \tensor \pi_1^*L \tensor \pi_2^*L^d)$.

Now consider a hypersurface $X$ in a $g$-dimensional abelian variety $A$ and let $L=\mathcal{O}_A(X)$.
\begin{lemma}
	If $L$ is ample then the surjectivity of the multiplication map
	\[
	H^0(A,L) \tensor H^0(A,L^{g-1}) \to H^0(A,L^g)
	\]
	implies the ITT.
	\begin{proof}
		Using the fact that the cotangent bundle of an abelian variety is trivial and that for an ample line bundle $L$ we have $H^i(A,L) = 0$ for $i>0$, it is easy to check that in each instance in Green's proof where $L$ is required to be sufficiently ample, ampleness is enough.
	\end{proof}
\end{lemma}
However, $L$ simply being ample is not sufficient to ensure the vanishing of 
$H^1(A \times A, \mathcal{I}_{\Delta/A} \tensor \pi_1^*L \tensor \pi_2^*L^{g-1})$. We will study more generally the surjectivity of the multipliction maps
\[
	\mu_n\colon H^0(A,L) \tensor H^0(A,L^n) \to H^0(A,L^{n+1})
\]
for $n \in \NN$.

\section{Surjectivity of multiplication maps}

A concept related to the surjectivity of the multiplication maps $\mu_n$ is projective normality. It can be 
defined for any projective variety and thus in particular for abelian varieties. We use the definitions given in \cite{MR2062673}.
\begin{dfn}
	A projective variety $Y \subset \PP^N$ is called \textit{projectively normal} in $\PP^N$ if its homogeneous 	
	coordinate ring is an integrally closed domain. A line bundle $L \rightarrow Y$ is called 
	\textit{normally generated} if it is very ample and $Y$ is projectively normal under the associated projective 
	embedding.
\end{dfn}
We can relate projective normality and the surjectivity of $\mu_n$. An ample line bundle $L$ on a projective variety $Y$ is normally generated if and only if the mutliplication map $H^0(Y,L) \tensor H^0(Y,L^n) \rightarrow H^0(Y,L^{n+1})$ is surjective for every $n \geq 1$ (see \cite[Lemma 7.3.2]{MR2062673}). 

This works for any projective variety but for abelian varieties surjectivity of $\mu_n$ implies surjectivity of $\mu_m$ for all $m \geq n$ (see for example \cite{MR1974682}). In particular we have that surjectivity of $\mu_1$ is equivalent to projective normality and that projective normality implies the ITT.

It is well known that a line bundle $L=M^n$ with $n \geq 3$ is normally generated and a line bundle $L=M^2$ is normally generated if and only if some additional assumption on the basepoints of $L$ holds. If we only care about surjectivity of $\mu_{g-1}$ the assumption on basepoints can be dropped at least when $g \geq 3$. 

Recall that for a line bundle $L$ on an abelian variety $A = V/\Lambda$ the first Chern class $c_1(L)$ defines a hermitian form on $V$ whose imaginary part $E$ is integer valued on $\Lambda$. The elementary divisor theorem ensures that there is a basis of $\Lambda$ such that $E$ is given by the matrix $\begin{pmatrix}
	0 	& D \\
	-D 	& 0
\end{pmatrix}$
where $D=\mathrm{diag}(d_1,\ldots,d_g)$ with integers $d_i \geq 0$ satisfying $d_i|d_{i+1}$ for $i=1,\ldots g-1$. The vector $(d_1,\ldots,d_g)$ is called the type of the line bundle $L$. If $L$ is ample we have $h^0(A,L) = d_1\cdots d_g$. Since $c_1(L^n)=nc_1(L)$ for any $n \in \NN$ by the above discussion the ITT holds for any smooth divisor in the linear system of a line bundle of type $(d_1, \ldots, d_g)$ with $d_1 \geq 2$.

The question remains what happens for primitive line bundles, i.e. line bundles of type $(1,d_2,\ldots,d_g)$. Note that by the Riemann-Roch theorem we have $h^0(A,L) = (L^g)/g!$ so any numerical condition can be equivalently expressed in terms of the number of sections of $L$ or the top intersection number.

For polarized abelian surfaces projective normality is fully understood. By \cite{Laz} and \cite{MR2076454}, if $(A,L)$ is a polarized abelian surface with $L$ very ample and of type $(1,d)$, then $L$ defines a projectively normal embedding if and only if $d>6$. Lazarsfeld's paper is hard to find but \cite{MR3656291} summarizes the main points. We already know that the ITT fails exactly on the locus of hyperelliptic curves with genus greater than $2$. By \cite[Theorem 2.8]{MR3968899} for any smooth hyperelliptic curve $C$ embedded in an abelian surface $A$ the genus $g(C)$ is $2,3,4$ or $5$ and $A$ is polarized of type $(1,g(C)-1)$. So smooth hypersurfaces of type $(1,5)$ and $(1,6)$ do not define projectively normal embeddings but satisfy the ITT. In the case that $g(C) = 2$ the ITT holds. By the above, $A$ is then principally polarized. The multiplication map $\mu$ cannot be surjective for purely dimensional reasons. This is however not a contradiction, as failure of $\mu$ to be surjective does not imply failure of the ITT.

For higher dimensional polarized abelian varieties Hwang and To show that a general polarized $g$-dimensional abelian variety with 
$h^0(A,L) \geq \frac{8^g}{2}\cdot \frac{g^g}{g!}$ is projectively normal (see \cite{MR2964474}). If we only want $\mu_{g-1}$ to be surjective we can in fact generalize the methods used in their proof to obtain a better bound.

We would prefer a different more explicit condition that we can check. Recall that an abelian variety $A$ is called \textit{simple} if the only abelian subvarieties are $\{0\}$ and $A$ itself. Iyer proves the following theorem.
\begin{thm}[\cite{MR1974682}]
	Let $L$ be an ample line bundle on a $g$-dimensional simple abelian variety $A$. If $h^0(A,L) > 2^g \cdot g!$, then $L$ gives a projectivelynormal embedding.
\end{thm}
Asymptotically this bound is worse than the one in \cite{MR2964474}. It does give a better bound up to $g = 23$. However the main reason we prefer this is that simplicity is a more conrete condition to work with.

This already gives us a sufficient condition for the ITT to hold but we can relax it. We cannot remove the condition that $A$ be simple by making the numerical condition on the global sections of $L$ stronger, even if we only try to prove surjectivity of $\mu_{g-1}$. In fact an analogous statement for any abelian variety cannot hold. Consider the abelian variety $X=C \times A$ where $(C,\mathcal{O}_C(2p))$ is a $(2)$-polarized elliptic curve and $(A,L)$ is a $(g-1)$-dimensional polarized abelian variety with polarization of type $(d_2,\ldots,d_g)$ where all $d_i$ are odd, e.g. a third power of a principal polarization. Now $X$ carries 
the product polarization $\mathcal{O}_C(2p)\boxtimes L$ which must be primitive because $\gcd(2,3) = 1$ but it cannot be normally generated as the restriction to $C$ is only basepoint free but not very ample. One would expect that for each abelian subvariety $B$ a numerical condition on the sections of the restriction $L|_B$ implying projective normality can be derived. Indeed, in the case of abelian threefolds Lozovanu proves the following theorem.
\begin{thm}[\cite{2018arXiv180308780L}]
		Let $(A,L)$ be a polarized abelian threefold such that $h^0(A,L) > 78$. Assume the following conditions:
		\begin{enumerate}
			\item For any abelian surface $S \subseteq A$ one has $h^0(S,L|_S) > 4$.
			\item For any elliptic curve $E \subseteq A$ one has $h^0(E,L|_E) > 4$.
		\end{enumerate}
		Then $L$ gives a projectively normal embedding of $A$.
\end{thm}
Note that he actually proves a more general result about $(A,L)$ satisfying the property $(N_p)$, however $(N_0)$ corresponds to projective normality.

\section{Proof of the main theorem}
	\label{section_main_thm}

	At the heart of the proof of Theorem \ref{thm_main_2} is the following lemma. This is the only place where $A$ needs to be simple.
	\begin{comment}
			\begin{lemma}[\cite{MR1974682}]
			\label{lemma_iyer}
			Let $D$ be an ample divisor on a $g$-dimensional simple abelian variety $A$. Let $G$ be a finite subgroup of $A$ contained in $D$. Then 
			$|G| \leq D^g$.
			\end{lemma}
			This can equivalently be written in a more geometric way.
	\end{comment}

	\begin{lemma}[\cite{MR1974682}]
		\label{prop_iyer}
		Let $L$ be an ample line bundle on a $g$-dimensional simple abelian variety $A$. Let $G$ be a finite subgroup with $|G| > h^0(A,L)\cdot g!$. Then the image of $G$ under the rational map $\phi_L\colon A \to \PP(H^0(A,L))^\dual$ generates $\PP(H^0(A,L))^\dual$.
	\end{lemma}
	Before going into the proof of Theorem \ref{thm_main_2} we recall some basic facts about polarized abelian varieties. Let $A = V/\Lambda$ be an abelian variety. A line bundle $L$ on $A$ induces a morphism 
	\begin{align*}
		\psi_L\colon & A \to \Pic^0(A) \\
					 & a \mapsto t_a^*L \tensor L^{-1}.
	\end{align*}
	Denote its kernel by $K(L)$. If $L$ is ample $\psi_L$ is an isogeny so that $K(L)$ is finite. A decomposition $\Lambda = \Lambda_1 \oplus \Lambda_2$ is a \textit{decomposition for $L$} if $\Lambda_1$ and $\Lambda_2$ are maximally isotropic with respect to the alternating form $\mathrm{Im} \ c_1(L)$. A decomposition $V = V_1 \oplus V_2$ is called a \textit{decomposition for $L$} if the induced decomposition $(V_1 \intersection \Lambda) \oplus (V_2 \intersection \Lambda)$ is a decomposition for $L$. Such a decomposition induces a decomposition $K(L) = K(L)_1 \oplus K(L)_2$.
	
	In the following let $(B,M)$ be a principally polarized abelian variety with $\theta \in H^0(B,M)$ the unique (up to a scalar) section. Write $B=V/\Lambda$ and let $\Lambda = \Lambda_1 \oplus \Lambda_2$ be a decomposition for $M$. Fix $n \in \NN$. 
	There is a natural action on $H^0(B,M^n)$ by the theta group $\mathcal{G}(M^n) = \{(b,\varphi) \ | \ b \in K(M^n), \ \varphi\colon t_b^*M^n \overset{\cong}{\to} M^n\}$. We can choose compatible isomorphisms $\varphi_b\colon t_b^*M^n \to M^n$ for $b \in K(M^m)_1$ so that for any $b,b' \in K(M^n)_1$ we have $\varphi_b(t_{b'}^*\varphi_{b'}) = \varphi_{b'}(t_b^*\varphi_b)$. That means that the action of $\mathcal{G}(M^n)$ induces an action of $K(M^n)_1$.
	For our purpose we want to find a section $\widetilde{\theta} \in H^0(B,M^n)$ that is invariant under this action. 
	Consider the isogeny
	\begin{align*}
	\varphi\colon B \to B'= B/K(M^n)_1
	\end{align*}
	and let $M'$ be a line bundle on $B'$ such that $\varphi^*M' = M^n$. Since $M'$ is a principal polarization there is a unique (again up to a scalar) section $\theta' \in H^0(B',M')$. We can take $\widetilde{\theta} = \varphi^*\theta'$ since clearly for any $\lambda \in K(M^n)_1$ we have $t_\lambda^*\widetilde{\theta} = t_\lambda^* \varphi^* \theta' = \varphi^* t_\lambda^* \theta' = \varphi^* \theta' = \widetilde{\theta}$ for any $\lambda \in K(M^n)_1$.
	Abusing notation a little we will also write $\theta$ and $\widetilde{\theta}$ for the associated theta divisors.
	
	Using the Theorem of the Square we see that for any $b \in B$
	\begin{align*}
		t^*_{nb}M \tensor t^*_{-b}M^n 	& \cong t_b^* M^n \tensor M^{-n+1} \tensor t_{-b}^*M^n \\
										& \cong M^n \tensor M^n \tensor M^{-n+1} \\
										& \cong M^{n+1}		
	\end{align*}
	so the divisor $t^*_{nb}\theta + t^*_{-b}\widetilde{\theta}$ is an element of the linear system $|(n+1)\theta|$ thus we have a morphism
	\begin{align*}
	\phi\colon 	& B \rightarrow  |(n+1)\theta| \\
	& b \mapsto t^*_{nb}\theta + t^*_{-b}\widetilde{\theta}.
	\end{align*}
	The following proposition is a generalization of a result by Wirtinger that can be found in \cite[p. 335]{MR0379510}.
	\begin{prop}
		\label{prop_1}
		For any $n \in \NN$ there is a nondegenerate bilinear form 
		$\eta\colon H^0(B,M^{n+1}) \tensor H^0(B,M^{n+1}) \rightarrow \CC$ inducing the isomorphism
		\[
		\eta'\colon \PP(H^0(B,M^{n+1})^\dual) \overset{\cong}{\longrightarrow} \PP(H^0(B,M^{n+1})) = |(n+1)\theta|
		\]
		such that
		\begin{equation*}
		\begin{tikzcd} 
		\												& \PP(H^0(B,M^{n+1})^\dual) \arrow{dd}{\eta'} \\
		B \arrow{ur}[swap]{\phi_{M^{n+1}}} \arrow{dr}{\phi}  & \ \\
		\												& \ |(n+1)\theta|
		\end{tikzcd}
		\end{equation*}
		commutes.
		\begin{proof}
			Consider the morphism
			\begin{align*}
			s\colon B \times B & \to B \times B \\
			(x,y)	& \mapsto (x+ny,x-y).
			\end{align*}
			We now have an isomorphism
			\[
			s^*(p_1^*M \tensor p_2^*M^n) \cong p_1^*M^{n+1} \tensor p_2^*M^{n(n+1)}
			\]
			To see this using the Appel-Humbert Theorem it suffices to compare the first Chern class and the semicharacters of both line bundles 
			(see \cite[Lemma 7.1.1]{MR2062673} for the case $n=1$). For any $m \in \NN$ and $\alpha \in K(M^m)_1$ we will write
			$\theta_\alpha^m = \varphi_\alpha(t_\alpha^*\theta^m)$ with $\varphi_\alpha \colon t_\alpha^*M^m \to M^m$ the compatibly chosen isomorphisms from before so that $\{\theta_\alpha^m \ | \ \alpha \in K(M^m)_1 \}$ defines a basis for 
			$H^0(B,M^m)$. Now we can write
			\begin{align}
			\label{eq_1}
			s^*(p_1^*\theta \tensor p_2^*\tilde{\theta}) = \sum_{\substack{\alpha \in K(M^{n+1})_1 \\ \beta \in K(M^{n(n+1)})_1}}
			{c_{\alpha\beta} p_1^*\theta_\alpha^{n+1} \tensor p_2^*\theta_\beta^{n(n+1)}}
			\end{align}
			We want to obtain dependencies between the coefficients $c_{\alpha\beta}$ to see that they are determined by a square matrix which we will use to define $\eta$. Consider the pullback of equation (\ref{eq_1}) by $t_{(0,-\gamma)}$ with $\gamma \in K(M^n)_1$. 
			On the left hand side, since
			\begin{align*}
			s(t_{(0,-\gamma)}(x,y)) &= s(x,y-\gamma) \\
			&= (x+n(y-\gamma),x-(y-\gamma))\\
			&= (x+ny-n\gamma,x-y+\gamma) \\
			&=t_{(0,\gamma)}(s(x,y))
			\end{align*}
			we get
			\begin{align*}
			t_{(0,-\gamma)}^*s^*(p_1^*\theta \tensor p_2^*\tilde{\theta}) 
			& = s^*t_{(0,\gamma)}^*(p_1^*\theta \tensor p_2^*\tilde{\theta}) \\
			& = s^*(p_1^*\theta \tensor p_2^*t_\gamma^*\tilde{\theta}) \\
			& = s^*(p_1^*\theta \tensor p_2^*\widetilde{\theta}).
			\end{align*}
			Here, we obtain the last line because we chose $\widetilde{\theta}$ such that it is invariant under translation by $\gamma \in K(M^n)_1$.
			On the right hand side we have
			\begin{align*}
			&	t_{(0,-\gamma)}^*\left(\sum_{\substack{\alpha \in K(M^{n+1})_1 \\ \beta \in K(M^{n(n+1)})_1}}
			{c_{\alpha\beta} p_1^*\theta_\alpha^{n+1} \tensor p_2^*\theta_\beta^{n(n+1)}} \right) \\
			= &	\sum_{\substack{\alpha \in K(M^{n+1})_1 \\ \beta \in K(M^{n(n+1)})_1}}
			{c_{\alpha\beta} p_1^*\theta_\alpha^{n+1} \tensor p_2^*t_{-\gamma}^*\theta_\beta^{n(n+1)}}
			\end{align*}
			The pullbacks on the right hand side permute the basis elements, comparing coefficients gives $c_{\alpha\beta} = c_{\alpha,\beta-\gamma}$. 
		
			Now because $\gcd(n,n+1)=1$, the exact sequence
			\[
				0 \to K(M^n)_1 \to K(M^{n(n+1)})_1 \to K(M^{n+1})_1 \to 0
			\]
			splits and thus $K(M^{n(n+1)})_1 \cong K(M^n)_1 \oplus K(M^{n+1})_1$. Therefore for any $\beta \in K(M^{n(n+1)})_1$ there is exactly one $\gamma \in K(M^n)_1$ such that $\beta-\gamma \in K(M^{n+1})_1$, namely $\gamma$ is the $n$-torsion part of $\beta$. Ultimately this means that we can choose representatives $\alpha,\beta \in K(M^{n+1})_1$ so that the matrix $(c_{\alpha\beta})$ is determined by $\alpha, \beta \in K(M^{n+1})_1$.
			
			\begin{comment}
				Translation by an element $b \in K(M^n)$ acts on $H^0(B,M^{n+1})$ and since $K(M^n) \subset K(M^{n(n+1)})$ there is an action of 
				$\Delta(B_{n+1}) = \{(b,b) \ | \ b \in B_{n+1}\}$ on $H^0(B,M^{n+1}) \tensor H^0(B,M^{n(n+1)})$. The element $s^*(p_1^*\theta \tensor p_2^*\widetilde{\theta})$ is invariant under this action and since the action on $H^0(B,M^{n+1})$ is irreducible it cannot lie in a proper subspace $W_1 \tensor W_2$ of $H^0(B,M^{n+1}) \tensor H^0(B,M^{n(n+1)})$. But if $\det(c_{\alpha\beta})$ were zero we could find bases $(v_1,\ldots,v_{(n+1)^g})$ and $(v'_1,\ldots,v'_{(n(n+1))^g})$ of $H^0(B,M^{n+1})$ and $H^0(B,M^{n(n+1)})$ respectively such that $s^*(p_1^*\theta \tensor p_2^*\tilde{\theta}) = \sum_{i,j}^r \alpha_{ij} p_1^*v_i\tensor p_2^*v'_j$ where the rank $r$ of $(c_{\alpha\beta})$ is strictly smaller than $(n+1)^g$ and thus $s^*(p_1^*\theta \tensor p_2^*\tilde{\theta}) \in \langle v_1,\ldots,v_r \rangle \tensor \langle v'_1,\ldots,v'_r \rangle$.
			\end{comment}
			We still need to show that $\det(c_{\alpha\beta}) \neq 0$. If the determinant were zero, the element 
			$s^*(p_1^*\theta \tensor p_2^*\tilde{\theta})$ would be contained in a proper subspace $W_1 \tensor W_2$ with $W_1 \subsetneq H^0(B,M^{n+1})$ of $H^0(B,M^{n+1}) \tensor H^0(B,M^{n(n+1)})$. However, translation by an element $b \in K(M^n)$ acts on $H^0(B,M^{n+1})$ and since $K(M^{n+1}) \subset K(M^{n(n+1)})$ there is an action of $\Delta(B_{n+1}) = \{(b,b) \ | \ b \in B_{n+1}\}$ on $H^0(B,M^{n+1}) \tensor H^0(B,M^{n(n+1)})$. The element $s^*(p_1^*\theta \tensor p_2^*\widetilde{\theta})$ is invariant under this action and since the action on $H^0(B,M^{n+1})$ is irreducible it cannot lie in such a proper subspace. 
			We conclude that $\det(c_{\alpha\beta}) \neq 0$ so $\eta(\theta_\alpha^{n+1},\theta_\beta^{n+1}) \coloneqq c_{\alpha\beta}$ defines the desired form $\eta$.
			
			The equation (\ref{eq_1}) can be expressed as
			\[
				\theta(u+nv)\tilde{\theta}(u-v) = \sum_{\alpha,\beta \in K(M^{n+1})_1}
				{c_{\alpha\beta} \theta_\alpha^{n+1}(u) \theta_\beta^{n(n+1)}(v)} \text{ for any } u,v \in B.
			\]
			For each $v \in B$ this implies that $u$ is in the support of the divisor $t_{nv}^* \theta + t_{-v}^* \widetilde{\theta}$ if and only if it is a zero of $\sum{c_{\alpha\beta}\theta_\beta^{n(n+1)}(v)\theta_\alpha^{n+1}}$ which gives that $\phi(v)=\eta'(\phi_M(v))$.
		\end{proof}
	\end{prop}
	
	With this we can prove the following theorem.
	\begin{thm}
		Let $L$ be a line bundle on a simple abelian variety $A$ of dimension $g$.
		Fix $n \in \NN$. If $h^0(A,L)>\left(\frac{n+1}{n}\right)^g \cdot g!$ then the multiplication map
		\[
		\mu_n\colon H^0(A,L) \tensor H^0(A,L^n) \rightarrow H^0(A,L^{n+1})
		\]
		is surjective.
	\begin{proof}
			Choose a maximal isotropic subgroup with respect to the Weil form, say $H=K(L)_1$ and cosider the isogeny
			\[
				\pi\colon A \to B=A/H.
			\]
			There is a principal polarization $M$ on $B$ such that $\pi^*M = L$. The character group $\widehat{H} \coloneqq \Hom(H,\CC^*)$ is a subgroup of $\Pic^0(B)$ so a character $\alpha \in \widehat{H}$ corresponds to a degree $0$ line bundle on $B$ also denoted by $\alpha$. 
			We have a decomposition $\pi_*\mathcal{O}_A = \bigoplus_{\alpha \in \widehat{H}}\alpha$. This gives us
			\begin{align*}
			\pi_* L 	& = \pi_*(\mathcal{O}_A \tensor L) \\
			& = \pi_*(\mathcal{O}_A \tensor \pi^*M) \\
			& = \pi_* \mathcal{O}_A \tensor M \hspace{3cm} \text{(projection formula)} \\
			& = \bigoplus_{\alpha \in \widehat{H}} M \tensor \alpha.
			\end{align*}
			More generally, for any $m \in \NN$, $\pi_* L^m = \bigoplus_{\alpha \in \widehat{H}} M^m \tensor \alpha$. Consequently
			\[
			H^0(A,L^m) \cong \bigoplus_{\alpha \in \widehat{H}} H^0(B,M^m \tensor \alpha).
			\]
			for any $m \in \NN$.			
			However, given a power of $L$ we take the larger subgroup $K(L^n)_1$ and get a finer decomposition. We will do that specifically for the second factor of $\mu_n$. Analogously to before, let $G=K(L^n)_1$ and consider the isogeny
			\[
			\pi'\colon A \to B' = A/G.
			\]
			Once again $B'$ is principally polarized say with polarization $M'$ and $L^n = \pi^*M'$. With the same arguments as above we can decompose
			\[
			H^0(A,L^n) \cong \bigoplus_{\alpha \in \widehat{G}} H^0(B,M' \tensor \alpha).
			\]
			Due to our choices of subgroups $H=K(L)_1=nK(L^n)_1$ is a subgroup of 
			$G$ so that these decompositions are compatible. This gives $\widehat{H} = \psi_M(\pi(K(L)_2))$ and $\widehat{G} = \psi_{M'}(\pi'(K(L^n)_2))$.
			The following diagram summarizes the situation
			\begin{equation}
			\label{cd_9}
			\begin{tikzcd}
			A \arrow{r}{\pi} \arrow{d}{\psi_L} \arrow[rr, bend left, "\pi'"]  & B \arrow{r}{\varphi} \arrow{d}{\psi_M} & B' \arrow{d}{\psi_{M'}} \\
			\Pic^0(A) 								& \Pic^0(B) \arrow{l}{\pi^*} 				 & \Pic^0(B') \arrow{l}{\varphi^*}.
			\end{tikzcd}
			\end{equation}
			Note that the second square does not commute but that we have instead $\varphi^* \circ \psi_{M'} \circ \varphi = n \cdot \psi_M$.
			
			Now we can write our multiplication map as
			\[
			\mu_n\colon \bigoplus_{\alpha \in \widehat{H}, \beta \in \widehat{G}} H^0(B,M \tensor \alpha) \tensor H^0(B',M' \tensor \beta) \overset{1\tensor\varphi^*}{\to}  
			\bigoplus_{\gamma \in \widehat{H}} H^0(B,M^{n+1}  \tensor \gamma).
			\]
			We can decompose $\mu_n = \bigoplus_{\gamma \in \widehat{H}}\mu_{n,\gamma}$ with
			\[
			\mu_{n,\gamma}\colon \bigoplus_{\beta \in \widehat{G}} H^0(B,M\tensor \gamma \tensor \varphi^*\beta) \tensor H^0(B',M' \tensor \beta^{-1}) \to H^0(B,M^{n+1} \tensor \gamma).
			\]
			Now since $\psi_M$ is an isomorphism we can take $H' \coloneqq \psi_M^{-1}(\widehat{H}) = \pi(K(L)_2)$,  
			$G'\coloneqq \psi_{M'}^{-1}(\widehat{G}) = \pi'(K(L^n)_2)$ and $\widetilde{G} \coloneqq \varphi^{-1}(G') \intersection \pi(K(L^n)_2)$
			Taking $c \in H'$ such that $\gamma = \psi_M((n+1)c) = \psi_{M^{n+1}}(c)$ and writing out the definitions of $\psi_M$ and $\psi_{M'}$, we obtain
			\begin{comment}
				\[ 
				\mu_{n,\gamma}\colon \bigoplus_{b \in G'} H^0(M\tensor \psi_M((n+1)c+nb)) \tensor H^0(M' \tensor \psi_{M'}(-b')) \to H^0(M^{n+1} \tensor \psi_{M^{n+1}}(c)).
							\]
				Finally writing out the definitions of $\psi_M$ and $\psi_{M'}$, we obtain
			\end{comment}
			\[
			\mu_{n,\gamma}\colon \bigoplus_{\substack{b' \in G',  b \in \widetilde{G} \\ \varphi(b) = b'}} H^0(B,t_{(n+1)c+nb}^*M) \tensor H^0(B',t_{-b'}^*M') \to H^0(B,t_c^*M^{n+1}).	
			\]	
			The difference between this and the proof in \cite{MR1974682} is that we are now taking the sum over the much larger group $G'$. Let $\theta$ be the unique theta divisor of $|M|$ and $\widetilde{\theta} \in |M^{n+1}|$ the pullback along $\varphi$ of the unique theta divisor $\theta'$ in $|M'|$. We see that $\mu_{n,\gamma}$ is surjective if the linear system $|t_c^*M^{n+1}|$ is generated by divisors of the form $t_{(n+1)c+nb}^*\theta+t_{-b}^*\widetilde{\theta} = t_c^*( t_{n(c+b)}^*\theta+t_{-(c+b)}^*\widetilde{\theta})$ with $b \in \widetilde{G}$. By Proposition \ref{prop_1} it is thus surjective if the image of $\widetilde{G}$ under $\phi_c\coloneqq t_c^* \circ\phi$ generates $|t_c^*M^{n+1}|$ or equivalently if the image of $\widetilde{G}$ under $\phi$ generates $|M^{n+1}|$. 
			Now by assumption we have 
			$|\widetilde{G}| = h^0(A,L^n) = n^g\cdot h^0(A,L)>(n+1)^g\cdot g! = h^0(B,M^{n+1})\cdot g!$ and thus we can apply Proposition \ref{prop_iyer} to finish the proof.  		
		\end{proof}
	\end{thm}
	Setting $n = g-1$ and using the method discussed in Section \ref{section_green} we obtain Theorem \ref{thm_main} as a corollary.
	\begin{corollary}
		\label{corollary_ITT}
		Let $X$ be a smooth hypersurface in a $g$-dimensional simple abelian variety $A$. If 
		$h^0(A,\mathcal{O}_A(X))>\left(\frac{g}{g-1}\right)^g \cdot g!$, then the Infinitesimal Torelli Theorem holds for $X$.
	\end{corollary}
	For the case $g=2$ this is exactly the same as in \cite{MR1974682}. However for higher dimensions our result directly improves the bound. For $g=3$ for example, $h^0(A,\mathcal{O}_A(X))>20$ is a sufficient condition for a hypersurface on a simple abelian variety to satisfy the ITT as we have seen above whereas to show that it gives a projectively normal embedding we need $h^0(A,\mathcal{O}_A(X))>48$.
	\begin{corollary}
		Let $S\subset A$ be a smooth complex projective surface that embeds into its Albanese $A$ as a hypersurface. If $S$ has geometric genus $p_g > 22$ and $A$ is simple then the ITT holds for $S$.
	\begin{proof}
		Consider the exact sequence
		\[
			0 \to \mathcal{O}_A \to \mathcal{O}_A(S) \to \mathcal{O}_S(S) \to 0
		\]
		By adjunction we have $\omega_S \cong \mathcal{O}_S(S)$ so taking cohomology and comparing dimensions gives
		\[
			h^0(\mathcal{O}_A(S)) = p_g + 1 - 3 > 20
		\]
		so we can apply Corollary \ref{corollary_ITT}.
	\end{proof}		
	\end{corollary}
	In \cite{MR929539} the ITT is proved for surfaces of irregularity greater than or equal to $5$ under the assumption that $\Omega_S^1$ is globally generated and that some other conditions hold. In our case $S$ has irregularity $3$ and we have to assume that the Albanese morphism $a:S \to A$ is an embedding which does in fact imply that $\Omega_S^1$ is globally generated. An interesting question would be if our approach can still be used to show the ITT when $\Omega_S^1$ is globally generated but $a$ is not an embedding. 
\newpage
\bibliography{The_Infinitesimal_Torelli_Theorem_for_hypersurfaces_in_abelian_varieties_-_Patrick_Bloss}{}

\providecommand{\bysame}{\leavevmode\hbox to3em{\hrulefill}\thinspace}
\providecommand{\MR}{\relax\ifhmode\unskip\space\fi MR }
% \MRhref is called by the amsart/book/proc definition of \MR.
\providecommand{\MRhref}[2]{%
  \href{http://www.ams.org/mathscinet-getitem?mr=#1}{#2}
}
\providecommand{\href}[2]{#2}
\begin{thebibliography}{{Loz}18}

\bibitem[Ago17]{MR3656291}
Daniele Agostini, \emph{A note on homogeneous ideals of abelian surfaces},
  Bull. Lond. Math. Soc. \textbf{49} (2017), no.~2, 220--225. \MR{3656291}

\bibitem[BL04]{MR2062673}
Christina Birkenhake and Herbert Lange, \emph{Complex abelian varieties},
  second ed., Grundlehren der Mathematischen Wissenschaften [Fundamental
  Principles of Mathematical Sciences], vol. 302, Springer-Verlag, Berlin,
  2004. \MR{2062673}

\bibitem[BO19]{MR3968899}
Pawe\l~\hspace{-6pt} Bor\'{o}wka and Angela Ortega, \emph{Hyperelliptic curves
  on {$(1,4)$}-polarised abelian surfaces}, Math. Z. \textbf{292} (2019),
  no.~1-2, 193--209. \MR{3968899}

\bibitem[Gar04]{MR2076454}
Luis~Fuentes Garc\'ia, \emph{Projective normality of abelian surfaces of type
  {$(1,2d)$}}, Manuscripta Math. \textbf{114} (2004), no.~3, 385--390.
  \MR{2076454}

\bibitem[GZ08]{MR2413341}
Umberto Garra and Francesco Zucconi, \emph{Very ampleness and the infinitesimal
  {T}orelli problem}, Math. Z. \textbf{260} (2008), no.~1, 31--46. \MR{2413341}

\bibitem[HT11]{MR2964474}
Jun-Muk Hwang and Wing-Keung To, \emph{Buser-{S}arnak invariant and projective
  normality of abelian varieties}, Complex and differential geometry, Springer
  Proc. Math., vol.~8, Springer, Heidelberg, 2011, pp.~157--170. \MR{2964474}

\bibitem[Iye03]{MR1974682}
Jaya~N. Iyer, \emph{Projective normality of abelian varieties}, Trans. Amer.
  Math. Soc. \textbf{355} (2003), no.~8, 3209--3216. \MR{1974682}

\bibitem[Laz90]{Laz}
Robert Lazarsfeld, \emph{Projectivit\'e normale des surfaces ab\'eliennes},
  (Redige par O. Debarre), Preprint No. 14, Europroj CIMPA (1990).

\bibitem[{Loz}18]{2018arXiv180308780L}
Victor {Lozovanu}, \emph{{Singular divisors and syzygies of polarized abelian
  threefolds}}, arXiv e-prints (2018), arXiv:1803.08780.

\bibitem[Mum74]{MR0379510}
David Mumford, \emph{Prym varieties. {I}}, Contributions to analysis (a
  collection of papers dedicated to {L}ipman {B}ers), 1974, pp.~325--350.
  \MR{0379510}

\bibitem[Rei88]{MR929539}
Igor Reider, \emph{On the infinitesimal {T}orelli theorem for certain irregular
  surfaces of general type}, Math. Ann. \textbf{280} (1988), no.~2, 285--302.
  \MR{929539}

\end{thebibliography}
\bibliographystyle{amsalpha}
%Comment out line below for appendix
%\end{document}
\newpage

\end{document}